\documentclass{article}
\usepackage{amssymb}
\usepackage{amsfonts}
\usepackage{amsmath}

\newtheorem{theorem}{Theorem}

\newtheorem{condition}[theorem]{Condition}

\begin{document}

\title{Survival of branching processes in random environments}
\author{E. Dyakonova\thanks{%
Steklov Mathematical Institute, Gubkina str. 8,
119991, Moscow, Russia. Email address: elena@mi.ras.ru. Supported in
part by the RFBR, Grant 11-01-00139.}, V. Vatutin\thanks{%
Steklov Mathematical Institute, Gubkina str. 8,
119991, Moscow, Russia. Email address: vatutin@mi.ras.ru. Supported
in part by the RFBR, Grant 11-01-00139.}, and S. Sagitov \thanks{%
Chalmers University and University of Gothenburg, 412 96 Gothenburg, Sweden. Email address:
serik@chalmers.se. Supported by the Swedish Research Council grant
621-2010-5623.}}
\maketitle

\begin{abstract}
This review paper presents the known results on the asymptotics of the
survival probability and limit theorems conditioned on survival of critical
and subcritical branching processes in IID random environments. The key
assumptions of the family of population models in question are:
non-overlapping generations, independent reproduction of particles within a
generation, independent reproduction laws between generations. This is a
biologically important generalization of the time inhomogeneous branching
processes. The assumption of IID (independent and identically distributed)
random environments reflects uncertainty in the future (as well as
historical) reproduction regimes in actual populations.

This review focusses on a particular range of questions of prime interest
for the authors. The reader should be aware of the fact that there are many
very interesting papers covering other issues on branching processes in
varying and random environments which are not mentioned here.
\end{abstract}

\section{Introduction}

There has been a constantly growing interest in branching processes in
random environment (BPREs) following the pioneering publications by Smith
and Wilkinsonin \cite{SmWil} in 1969 and Athreya and Karlin \cite{AK71}, 
\cite{AK71a} in 1971 (see \cite{VZ} for a list of references up to 1985 and 
\cite{af80}--\cite{Koz}, \cite{BGK}--\cite{de}, \cite{Dyak9}--\cite{DyakGVa}%
, \cite{FV}--\cite{GuLi01}, \cite{K99}, \cite{Liu}, \cite{DH}, and \cite{V01}%
--\cite{VW09} for some more recent results). A BPRE is a stochastic
population model where individuals constituting a generation reproduce
independently according to a common offspring distribution. The random
environment feature of the model means that the reproduction laws randomly
change from one generation to the other. In this review we summarize some
recent foundings on critical and subcritical discrete time BPREs under a key
assumption that the reproduction laws of particles of different generations
are independent and identically distributed (IID). The results reported here
are mainly due to the joined efforts of the Russian-Germain team of
researchers consisting of V. Afanasyev, C. Boinghoff, E. Dyakonova, J.
Geiger, K. Fleischmann, G. Kersting, V. Vatutin, \ V. Wachtel and others
(see also \cite{BGK} for an earlier review).

In the study of stochastic processes in random environment it is crucial to
recognize the difference between quenched and annealed approaches. Under the 
\textit{quenched approach} characteristics of a BPRE such as the survival
probability at time~$n$ are treated as random variables, where the source of
randomness is due to uncertainty in possible realizations of the environment
(see, for example,~ \cite{AK71}, \cite{AK71a}, \cite{AN}, \cite{Ta}, \cite%
{VDfr} and the relevant literature in~\cite{VZ}). While the \textit{annealed
approach} studies the mean values of the mentioned characteristics as a
result of averaging over possible realizations of the environment (see \cite%
{Af}, \cite{Af97}, \cite{Af3}, \cite{BV}, \cite{FV}, \cite{K76}, \cite{K99}, 
\cite{V01}, \cite{VD97} and the bibliography in~\cite{VZ}). The annealed
approach can be viewed as a summary of the more detailed quenched analysis.
One has to be aware that in biological applications the predictions by the
annealed approach might be misleading: in reality the environmental history,
being uncertain, is unique.

The remainder of the paper is organized as follows. Section \ref{Sclass}
introduces an important classification of BPREs into three major classes:
supercritical, critical, and subcritical branching processes. Here we
discuss the key assumptions imposed on the reproduction laws in the critical
and subcritical cases. Section \ref{Saa} contains asymptotic results
obtained using the annealed approach. In Section \ref{Smu} we summarize the
results for the multi-type BPREs where individuals within the same
generation are allowed to have different reproduction laws. The quenched
results for the critical BPREs are presented in Section \ref{Squ}. The last
Section 6 is devoted to the reduced critical BPREs.

Observe that in the setting of random environment we use different fonts for
the probability and expectation operators depending on whether we condition
on the environment ($\mathcal{P},\mathcal{E}$) or not ($\mathbb{P},\mathbb{E}
$). Yet another font ($\mathbf{P},\mathbf{E}$) is used when we introduce an
auxiliary probability measure in the subcritical case (compare for example
Conditions \ref{B2} and \ref{A3ws}). Constants denoted by $c$ may take
different values in different formulae (the same reservation holds for $%
c_1,c_2$).

\section{\protect\bigskip Classification of BPREs}

\label{Sclass}

The basic classification of the branching processes in a constant
environment recognizes three major reproduction regimes: subcritical,
critical and supercritical. It is based on the mean value $m$ of the
offspring number and these three classes are defined by the following
relations: $m<1$, $m=1$, and $m>1$. Adding the feature of random environment
makes the classification issue less straightforward due to random
fluctuations of the consecutive offspring mean values.

\subsection{BPREs and associated random walks}

A branching process in a \textit{varying} environment is most conveniently
described in terms of a sequence of probability generating functions (PGFs) 
\begin{equation*}
f_{n}(s)=\sum_{i=0}^{\infty }q_{i}^{\left( n\right) }s^{i},
\end{equation*}%
where $q_{i}^{\left( n\right) }$ stands for the probability that a particle
from generation $n-1$ contributes to the next generation $n$ by producing
exactly $i$ offspring (independently of other particles). Time inhomogeneous
Markov chain $\{Z_{n}\}$ describing the fluctuations in population sizes is
then characterized by the iterations of the consecutive reproduction PGFs 
\begin{equation}
\mathbb{E}\left[ s^{Z_{n}}\right] =f_{1}(f_{2}(\ldots (f_{n}(s))\ldots )),
\label{ve}
\end{equation}%
if we assume that $Z_{0}=1$. Taking derivatives in the previous relation
gives an expression for the expected population size 
\begin{equation}
\mathbb{E}\left[ Z_{n}\right] =m_{1}\ldots m_{n} ,  \label{mean}
\end{equation}%
where $m_{n}=f_{n}^{\prime }(1)$ stands for the mean offspring number.

An important example is the linear-fractional case 
\begin{equation*}
f_{n}(s)=r_{n}+(1-r_{n})\frac{t_{n}s}{1-(1-t_{n})s}
\end{equation*}%
fully characterized by a sequence of pairs of parameters $(r_{n},t_{n})\in
\lbrack 0,1)\times (0,1]$. The iteration of linear-fractional PGFs is again
linear-fractional 
\begin{equation*}
f_{1}(f_{2}(\ldots (f_{n}(s))\ldots ))=r^{(n)}+(1-r^{(n)})\frac{t^{(n)}s}{%
1-(1-t^{(n)})s}
\end{equation*}%
with 
\begin{align*}
\frac{1-r^{(n)}}{t^{(n)}}& =m_{1}\ldots m_{n}, \\
\frac{1-t^{(n)}}{t^{(n)}}& =\frac{1-t_{n}}{t_{n}}+\frac{1-t_{n-1}}{t_{n-1}}%
m_{n}+\ldots +\frac{1-t_{1}}{t_{1}}m_{2}\ldots m_{n}.
\end{align*}%
It follows, 
\begin{equation*}
\frac{1}{t^{(n)}}=1+\frac{1-t_{n}}{1-r_{n}}m_{n}+\frac{1-t_{n-1}}{1-r_{n-1}}%
m_{n-1}m_{n}+\ldots +\frac{1-t_{1}}{1-r_{1}}m_{1}\ldots m_{n},
\end{equation*}%
and since $\mathbb{P}(Z_{n}=0)=r^{(n)}$, 
\begin{equation*}
\frac{1}{\mathbb{P}(Z_{n}>0)}=e^{-s_{n}}+\sum_{k=1}^{n}a_{k}e^{-s_{k-1}},
\end{equation*}%
with $a_{k}=\frac{1-t_{k}}{1-r_{k}}$ and $s_{k}=\ln m_{1}+\ldots +\ln m_{k}$%
. The obtained equality indicates that the asymptotics of the survival
probability $\mathbb{P}(Z_{n}>0)$ under certain assumptions is governed by
the minimal value of consecutive expectations of population sizes 
\begin{equation*}
\min (e^{s_{1}},\ldots ,e^{s_{n}})=e^{\min (s_{1},\ldots ,s_{n})}.
\end{equation*}%
Notice also that in the linear-fractional case $(Z_{n}-1\mid Z_{n}>0)$ is
geometrically distributed Geom$(t^{(n)})$.

Given a sequence of \textit{random} PGFs $(F_{1},F_{2},...)$ one can speak
of a BPRE $\{Z_n\}_{n=0}^\infty$. Throughout we will assume that

\begin{itemize}
\item the defining random PGFs in are IID so that all $F_n\overset{d}{=}F$
have the same marginal distribution,

\item the BPRE starts by a single particle, $Z_0=1$, unless it is clearly
stated otherwise.
\end{itemize}

For BPREs the relation \eqref{ve} transforms into 
\begin{equation*}
\mathbb{E}\left[ s^{Z_{n}}\mid F_{1}=f_{1},F_{2}=f_{2},\ldots \right]
=f_{1}(f_{2}(\ldots (f_{n}(s))\ldots )).
\end{equation*}%
Throughout we will denote the probability and expectation conditioned on the
environment by 
\begin{equation*}
\mathcal{P}(\cdot ):=\mathbb{P}\left( \cdot \mid F_{1},F_{2},\ldots \right)
,\quad \mathcal{E}\left[ \cdot \right] :=\mathbb{E}\left[ \cdot \mid
F_{1},F_{2},\ldots \right] .
\end{equation*}%
Observe that the quenched counterpart of the relation \eqref{mean} can be
viewed as a realization of a Markov chain $\mathcal{E}\left[ Z_{n}\right]
=e^{S_{n}}$ defined by the so called \emph{associated random walk} (ARW): 
\begin{equation*}
S_{0}=0,\ S_{n}=X_{1}+...+X_{n},\,\mbox{ where }X_{n}=\ln F_{n}^{\prime
}\left( 1\right) .
\end{equation*}%
Here the jump sizes are IID random variables with 
\begin{equation*}
X_{n}\overset{d}{=}X:=\ln F^{\prime }(1)=\ln \mathcal{E}\left[ \xi \right] ,
\end{equation*}%
where $\xi $ is a random variable representing the offspring number for a
single particle reproducing in a random environment.

The key representation $\mathcal{E}\left[ Z_{n}\right] =e^{S_{n}}$ leads to
the following extended classification of branching processes: a BPRE is
called (A) supercritical, (B) subcritical, (C) non-degenerate critical and (C%
$_{0}$) degenerate critical, if its ARW a.s. satisfies one of the following
conditions

(A) $\lim_{n\rightarrow \infty }S_{n}=+\infty $,

(B) $\lim_{n\rightarrow \infty }S_{n}=-\infty $,

(C) $\lim \sup_{n\rightarrow \infty }S_{n}=+\infty $ $\ $\ and $\lim
\inf_{n\rightarrow \infty }S_{n}=-\infty $,

(C$_0$) $S_{n}=0$ for all $n.$

\noindent This classification, due to \cite{4h}, is based on the crucial
fact (see for example, \cite{Fe}, ch.12, \S 2) that \textit{any} IID random
walk can be attributed to one of these four classes. For the processes in
which $\mathbb{E}\left[ X\right] $ exists this classification coincides with
the standard classification of branching processes: a BPRE is called
supercritical, critical or subcritical if $\mathbb{E}\left[ X\right] >0,$ $%
=0,$ or $<0,$ respectively (see, for instance, \cite{AN}). However, the
extended classification seems to be more natural and is justified by a
series of results published in \cite{4h}, \cite{DyakGVa}, \cite{VDfr}, \cite%
{VD}, \cite{VD2}, \cite{VDWie}. In particular, in view of the estimate 
\begin{equation*}
\mathcal{P}\left( Z_{n}>0\right) =\min_{1\leq k\leq n}\mathcal{P}%
(Z_{k}>0)\leq \min_{1\leq k\leq n}\mathcal{E}\left[ Z_{k}\right] =e^{\min
(S_{1},\ldots ,S_{n})}
\end{equation*}%
it follows that the probability of extinction of subcritical and
non-degenerate critical BPRE equals 1 a.s. Moreover, it was shown in \cite%
{4h}, \cite{DyakGVa}, \cite{VDfr}, \cite{VD}, \cite{VD2}, and \cite{VDWie}
that the asymptotic properties of the non-degenerate critical BPRE with $%
\mathbb{E}\left[ X\right] =0$ are quite similar to those in which $\mathbb{E}%
\left[ X\right] $ does not exist.

In what follows we mostly focus on classes (B) and (C) excluding the
degenerate critical class (C$_0$). Therefore, in the case (C) we will speak
of critical BPREs omitting the specification $^{\prime }$non-degenerate$%
^{\prime }$. Note that by excluding the class (C$_0$) we disregard the
critical Galton-Watson processes in constant environments.

\bigskip

\subsection{The Spitzer-Doney condition in the critical case}

We further specify the general condition (C) by assuming that the ARW
satisfies the classical \textit{Spitzer condition} with parameter $\rho $ : 
\begin{equation*}
\frac{1}{n}\sum_{k=1}^{n}\mathbb{P}(S_{k}>0)\rightarrow \rho \in (0,1)\quad 
\text{as}\quad n\rightarrow \infty .
\end{equation*}%
According to \cite{Do} the Spitzer condition is equivalent to the following 
\textit{Doney condition} with parameter $\rho $ 
\begin{equation*}
\mathbb{P}(S_{n}>0)\rightarrow \rho \in (0,1)\quad \text{as}\quad
n\rightarrow \infty .
\end{equation*}%
In this text we refer to this condition as the \textit{Doney-Spitzer
condition} expressing it in the form 
\begin{equation}
\rho =\lim_{n\rightarrow \infty }\frac{1}{n}\sum_{k=1}^{n}\mathbb{P}\left(
S_{k}>0\right) =\lim_{n\rightarrow \infty }\mathbb{P}\left( S_{n}>0\right)
\in (0,1).  \label{dpar}
\end{equation}%
It is well-known that any random walk satisfying \eqref{dpar} is of the
oscillating type (see, e.g., Section XII.7 in \cite{Fe}). We note that
condition \eqref{dpar} covers random walks satisfying Condition \ref{A1(a)}
(see below) as well as all non-degenerate symmetric random walks. In these
cases $\rho =1/2$.

\begin{condition}
\label{A1(a)} The distribution of $X$ has zero mean and variance $0<\sigma
^{2}<\infty $. It is non-lattice.
\end{condition}

Of course, one can consider the lattice case as well. However, this leads to
unnecessary complications in the statements of the respective results which
we prefer to avoid.

Other examples where the Doney-Spitzer condition fulfills are provided by
random walks in the domain of attraction of some stable law, see Condition %
\ref{B1} below. We will often refer to this weaker version of \eqref{dpar}
as it allows to relax some extra conditions imposed on the random
environment.

\begin{condition}
\label{B1} The distribution of $X$ is non-lattice and belongs without
centering to the domain of attraction of a stable law with index $\alpha \in
(0,2]$ and skewness parameter $\psi $, such that

\begin{itemize}
\item $|\psi |<1$, if $0<\alpha <1$,

\item $|\psi |\le1$, if $1<\alpha <2$,

\item $\psi=0$ if $\alpha =1$ or $\alpha =2$.
\end{itemize}
\end{condition}

Recall that the condition that the centering constants are zero is
redundant, if $0<\alpha <1$, and that a stable law from Condition \ref{B1}
has the characteristic function (see, e.g., Theorem~8.3.2 in~\cite{Bbook}) 
\begin{equation*}
\chi _{\alpha ,\psi }\left( t\right) =\exp \Big(-c|t|^{\,\alpha }\Big(%
1-i\psi \frac{t}{|t|}\tan \frac{\pi \alpha }{2}\Big)\Big),\ t\in \mathbb{R},
\label{std}
\end{equation*}%
where $c>0$ is a scaling parameter. This stable law has finite absolute
moments of all orders $r<\alpha $. One can check (see, e.g., Section 8.9.2
in~\cite{Bbook}) that Condition~\ref{B1} implies \eqref{dpar} with parameter 
\begin{equation}
\rho =\frac{1}{2}+\frac{1}{\pi \alpha }\arctan \left( \psi \tan \frac{\pi
\alpha }{2}\right) \in \,(0,1).  \label{par}
\end{equation}

Let%
\begin{equation*}
L_{n}:=\min \left( S_{1},...,S_{n}\right) ,\qquad M_{n}:=\max \left(
S_{1},...,S_{n}\right) ,
\end{equation*}%
and 
\begin{equation*}
\tau \left( n\right) :=\min \{i\in \lbrack 0,n]:\;S_{i}=\min \left(
0,L_{n}\right) \}
\end{equation*}%
be the left--most point at which the ARW attains its minimal value on the
time-interval $[0,n].$ It is known (see\ \cite[Ch.~IV, $\S \,$20]{S}) that
under the Doney-Spitzer condition 
\begin{equation}
n^{-1}\tau (n)\overset{d}{\rightarrow }\tau ,\qquad n\rightarrow \infty ,
\label{converarc}
\end{equation}%
where $\tau $~is a random variable having a Beta-distribution with
parameters\linebreak $(1-\rho ,\rho )$ and the symbol $\overset{d}{%
\rightarrow }$ stands for the convergence in distribution.

Let 
\begin{equation*}
D:=\sum_{k=1}^{\infty }k^{-1}\mathbb{P}\left( S_{k}=0\right),
\end{equation*}%
and 
\begin{align*}
&\gamma _{0}:=0,\quad \gamma _{j+1}:=\min (n>\gamma _{j}:S_{n}<S_{\gamma
_{j}}), \\
&\Gamma _{0}:=0,\quad \Gamma _{j+1}:=\min (n>\Gamma _{j}:S_{n}>S_{\Gamma
_{j}}),\,j\geq 0,
\end{align*}%
be the strict descending and strict ascending ladder epochs of the ARW. Next
we introduce two renewal functions 
%

\begin{equation}
U(x)=\left\{ 
\begin{array}{ccc}
1_{\{0<x\}}+\sum_{j=1}^{\infty }\mathbf{P}(S_{\Gamma _{j}}<x) & \hbox{if} & 
x>0, \\ 
e^{-D} & \hbox{if} & x=0, \\ 
0 & \hbox{if} & x<0,%
\end{array}%
\right.  \label{DevU}
\end{equation}%
and%
\begin{equation}
V(x)=\left\{ 
\begin{array}{ccc}
\sum_{j=0}^{\infty }\mathbf{P}(S_{\gamma _{j}}\geq -x) & \hbox{if} & x>0, \\ 
1 & \hbox{if} & x=0, \\ 
0 & \hbox{if} & x<0.%
\end{array}%
\right.  \label{DefV}
\end{equation}%
It is known (see \cite{VD2}, Lemma 1) that \eqref{dpar} implies 
\begin{eqnarray}
&&\mathbb{E}\left[ U(-X)1_{\{X<0\}}\right] =e^{-D},\quad  \notag \\
&&\mathbb{E}\left[ U(x-X)1_{\{X<x\}}\right] =U(x),\quad x>0,  \notag \\
&&\mathbb{E}\left[ V(x+X)\right] =V(x),\quad x\geq 0.  \label{PV}
\end{eqnarray}

\subsection{Refined classification in the subcritical case}

\label{Ss}

According to (B) the ARW in the subcritical case should have a clear trend
towards minus infinity. It turns out (see \cite{af80}, \cite{de}, and \cite%
{GKV}) that the behavior of a subcritical BPRE to a large extend is
determined by the speed of the negative drift quantified by the parameter $%
\beta $, when exists, such that 
\begin{equation}
\mathbb{E}\left[ Xe^{\beta X}\right] =0.  \label{beta}
\end{equation}%
We will distinguish among three different sub-cases:

\begin{itemize}
\item \emph{weakly} subcritical, if \eqref{beta} holds with some $0<\beta <1$%
, which implies $\mathbb{E}\left[ Xe^{X}\right] >0$,

\item \emph{intermediately} subcritical, if \eqref{beta} holds with $\beta
=1 $, so that $\mathbb{E}\left[ Xe^{X}\right] =0$,

\item \emph{strongly} subcritical, if $\mathbb{E}\left[ Xe^{X}\right] <0$.
\end{itemize}

Notice that due to monotonicity of the function $\mathbb{E}\left[ Xe^{\beta
X}\right] $ in all three sub-cases we have $\mathbb{E}\left[ X\right] <0$.
Clearly, the smaller is the positive value of $\beta $ the closer we get to
the critical case with $\mathbb{E}\left[ X\right] =0$. On the other hand,
the subcritical case in a constant environment corresponds to $X$ being a
negative constant implying the strong subcriticality.

Observe next that $\mathbb{E}\left[ e^{\beta S_{n}}\right] =\gamma ^{n}$,
where 
\begin{equation}
\gamma =\mathbb{E}[e^{\beta X}]\ .  \label{dgamm}
\end{equation}%
Let us introduce the following auxiliary measure $\mathbf{P}$ with
expectation $\mathbf{E}$. For any $n\in \mathbf{N}$ and any measurable,
bounded function $\varphi :\Delta ^{n}\times \mathbf{N}_{0}^{n+1}\rightarrow 
\mathbb{R}$, the measure $\mathbf{P}$ is given by 
\begin{equation*}
\mathbf{E}[\varphi (F_{1},\ldots ,F_{n},Z_{0},\ldots ,Z_{n})]\ :=\ \gamma
^{-n}\mathbb{E}\big[\varphi (F_{1},\ldots ,F_{n},Z_{0},\ldots
,Z_{n})e^{\beta S_{n}}\big].
\end{equation*}%
Notice that $\mathbf{E}[X]\ =\ 0$ making $S$ a recurrent random walk under $%
\mathbf{P}$. Clearly, in the critical case the measures $\mathbf{P}$ and $%
\mathbb{P}$ coincide since $\beta =0$.

\begin{condition}
\label{A2ws} The distribution of $X$ with respect to $\mathbf{P}$ is
non-lattice, has zero mean and belongs to the domain of attraction of a
stable law with index $\alpha \in (1,2]$.
\end{condition}

Under Condition \ref{A2ws} there exists an increasing sequence of positive
numbers $a_n$ regularly varying at infinity 
\begin{equation}
a_{n}\ =\ n^{1/\alpha }l_{n}  \label{defcoefal}
\end{equation}%
such that the scaled ARW $a_n^{-1}S_{nt}$ $\mathbf{P}$-weakly converges to a
strictly stable L\'{e}vy-process $L_t$ with parameter $\alpha$. Here and
elsewhere in the expressions like $S_{nt}$ the index $nt$ is understood as
its integer part.

\section{Annealed approach}

\label{Saa}

Different trajectories of the ARW represent different scenarios of
environmental history. The more favorable scenarios translate into higher
ARW trajectories. If an ARW trajectory stays below zero for a longer period
of time, the BPRE is doomed to die out. The annealed analysis of the
survival of BPRE gives a summary picture of the system behavior after
averaging over all successful ARW trajectories.

\subsection{\protect\bigskip Asymptotic behavior of the survival probability}

The first result on the asymptotic behavior of the survival probability for
the critical BPREs has appeared in the seminal paper by Kozlov \cite{K76}.
Kozlov has shown for BPREs satisfying Condition \ref{A1(a)} that under
certain additional assumptions 
\begin{equation*}
\qquad c_{1}n^{-1/2}\leq \mathbb{P}(Z_{n}>0)\leq c_{2}n^{-1/2},
\label{VD97(1.3)}
\end{equation*}%
where $0<c_1\le c_2<\infty.$ Additionally, he proved that in the
linear-fractional case 
\begin{equation}
\qquad \mathbb{P}(Z_{n}>0)\ \sim \ cn^{-1/2}\ ,\quad c>0,\quad n\rightarrow
\infty .  \label{asym}
\end{equation}%
These results were published in 1976. Only in 2000 Geiger and Kersting \cite%
{GK} were able to demonstrate the validity of (\ref{asym}) for arbitrary
BPREs satisfying Condition \ref{A1(a)}. To establish this result Geiger and
Kersting developed a new powerful method of proving conditional limit
theorems for BPREs based on a change of measures. This method, along with
the idea of splitting trajectories of the ARW at the point of global minimum
on the interval $\left[ 0,n\right] ,$ first suggested by Vatutin and
Dyakonova \cite{VD}, became one of the main tools for obtaining conditional
limit theorems not only for the critical BPREs but for the subcritical BPREs
as well.

\subsubsection{\protect\bigskip Critical case}

The main idea in finding the asymptotic behavior of the survival probability
of critical BPREs under annealed approach is to show that 
\begin{equation}
\mathbb{P}\left( Z_{n}>0\right) \ \sim \ \theta \,\mathbb{P}\left( \min
(S_{1},\ldots ,S_{n})\geq 0\right) ,\quad n\rightarrow \infty ,  \label{key}
\end{equation}%
where $\theta >0$. This is done by splitting trajectories of the ARW at the
point of global minimum on the interval $\left[ 0,n\right] $. For the
annealed approach the trajectory splitting method was first used by
Dyakonova, Geiger, and Vatutin \cite{DyakGVa} and then further developed by
Afanasyev, Geiger, Kersting, and Vatutin \cite{4h}.

Recall that $X=\ln \mathcal{E}\left[ \xi \right] $ is the logarithm of the
mean for the conditional on the environment offspring number distribution 
\begin{equation}
Q_{i}=\mathcal{P}(\xi =i),\ i=0,1,2,\ldots .  \label{Q}
\end{equation}%
In terms of the conditional probabilities $Q_{i}$ the random PGF $F$ can be
expressed as 
\begin{equation*}
F(s)=\mathcal{E}\left[ s^{\xi }\right] =\sum_{i=0}^{\infty }Q_{i}s^{i}.
\end{equation*}%
In \cite{4h} a higher moment assumption on the environment is given in terms
of the standardized truncated second moment of the offspring number $\xi $: 
\begin{equation}
\zeta (a)\ :=\ \frac{\mathcal{E}\left[ \xi ^{2}1_{\{\xi \geq a\}}\right] }{(%
\mathcal{E}\left[ \xi \right] )^{2}}.  \label{defdzet}
\end{equation}

\begin{condition}
\label{A2} For some positive $\varepsilon$ and $a$ 
\begin{equation*}
\mathbb{E}\,\left[ (\ln ^{+}\zeta (a))^{\frac{1}{\rho }+\epsilon }\right] \
<\ \infty \quad \text{and}\quad \mathbb{E}\,[V(X)(\ln ^{+}\zeta
(a))^{1+\epsilon }]\ <\ \infty ,
\end{equation*}%
where $\ln ^{+}x:=\ln \left( \max \left\{ x,1\right\} \right) $.
\end{condition}

Next we present a number of cases when this assumption is fulfilled.

\begin{enumerate}
\item If the offspring number has a bounded support, i.e. 
\begin{equation}
\mathbb{P}(\xi \leq a^{\ast })=1  \label{bounded}
\end{equation}%
for some $a^{\ast }$, then $\zeta (a)=0$ $\mathbb{P}$--a.s. for all $%
a>a^{\ast }$. Obviously, in this case Condition~\ref{A2} is valid. As a
particular example we mention here the binary splitting reproduction law
characterized by 
\begin{equation*}
\mathcal{P}(\xi=2)=1-\mathcal{P}(\xi=0)=e^{X}/2.
\end{equation*}

\item In view of relation (\ref{PV}), we have $\mathbb{E}\left[ V(X)\right]
=V(0)=1$. Therefore, Condition \ref{A2} is satisfied, if $\zeta (a)$ is a.s.
bounded from above for some $a\in \mathbf{N}_{0}$.

Observe that $\zeta (2)\leq 2\eta $, where 
\begin{equation}
\eta :=\frac{F^{\prime \prime }\left( 1\right) }{\left( F^{\prime }\left(
1\right) \right) ^{2}}=\frac{\mathcal{E}\left[ \xi (\xi -1)\right] }{\left( 
\mathcal{E}\left[ \xi \right] \right) ^{2}}  \label{ddefeta}
\end{equation}%
implying that $\zeta (2)$ is bounded from above if $\eta $ is bounded. This
is the case when \eqref{Q} is either a Poisson distribution with a random
mean (when $\eta =1$ a.s.) or a geometric distribution with a random mean
(when $\eta =2$ a.s.).

\item The renewal function $V(x)$ always satisfies $V(x)=O(x)$ as $%
x\rightarrow \infty $ and $V(x)=0$ for $x<0$. Therefore, as follows from H%
\"{o}lder's inequality, Condition \ref{A2} would be satisfied, provided 
\begin{equation*}
\mathbb{E}\left[ (X^{+})^{p}\right] \ <\ \infty \quad \text{and}\quad 
\mathbb{E}\left[ (\ln ^{+}\zeta (a))^{q}\right] <\infty \ 
\end{equation*}%
for some $p>1$ and $q>\max (\rho ^{-1},p(p-1)^{-1})$.
\end{enumerate}

Under Condition \ref{B1} we have $\rho \ \leq \ \alpha ^{-1}$ due to %
\eqref{par}. This observation leads to the following weaker version of
Condition \ref{A2}.

\begin{condition}
\label{B2} For some positive $\varepsilon$ and $a$ 
\begin{equation*}
\mathbb{E}\,\left[ (\ln ^{+}\zeta (a))^{\alpha +\epsilon }\right] \ <\
\infty \ .
\end{equation*}
\end{condition}

It was proved in \cite{4h} that if \eqref{dpar} holds together with
Condition \ref{A2} or if Condition \ref{B1} with $|\psi |<1$ and Condition %
\ref{B2} are valid, then there exists a positive finite number $\theta $
such that, as $n\rightarrow \infty ,$ 
\begin{equation}
\mathbb{P}\left( Z_{n}>0\right) \ \sim \ \theta \,\mathbb{P}\left( \min
(S_{1},\ldots ,S_{n})\geq 0\right) \ \sim \theta n^{-(1-\rho )}l_{n},
\label{ratio}
\end{equation}%
where $l_{1},l_{2},\ldots $ is a sequence varying slowly at infinity
(consult \cite{4h} for representations of $\theta $ and the sequence $l$).
Moreover, it was shown in \cite{4h} that for a trajectory of the ARW to give
an essential contribution to the annealed probability of survival during the
time interval $[0,n]$ this trajectory should pass the point of its global
minimum over this time interval very soon after its start at time 0.

The result \eqref{ratio} gives the asymptotics of the tail distribution for
the extinction time $T=\min \left\{ k:Z_{k}=0\right\} $ in view of the
equality 
\begin{equation*}
\mathbb{P}(Z_{n}>0)=\mathbb{P}(T>n).
\end{equation*}%
The study of the asymptotic behavior of $\mathbb{P}\left( T=n\right) $ for
the critical BPRE was initiated by Vatutin and Dyakonova \cite{VD97}. They
demonstrated that for the linear-fractional BPREs under \eqref{dpar},
Condition \ref{A1(a)} and some additional conditions one has 
\begin{equation}
\qquad \mathbb{P}\left( T=n\right) \ \sim \ cn^{-3/2}\ ,\quad n\rightarrow
\infty ,  \label{VD97(1.7)}
\end{equation}%
where $c\in (0,\infty ).$ Afterwards, B\"{o}inghoff, Dyakonova, Kersting,
and Vatutin \cite{BDKV} have proved (\ref{VD97(1.7)}) for the
linear-fractional BPREs under weaker conditions: \eqref{dpar}, Condition \ref%
{A1(a)}, Condition \ref{B2} with $\alpha=2$ plus the following Condition \ref%
{A0}.

\begin{condition}
\label{A0} There exists a constant $\chi \in (0,1/2)$ such that \textbf{\ }%
\begin{equation*}
\mathbb{P}(\chi \leq Q_{0}\leq 1-\chi ,\quad \eta \geq \chi )=1.
\label{Asss1}
\end{equation*}
\end{condition}


Vatutin and Wachtel \cite{VW09} for the geometric case with 
\begin{equation*}
F_{n}(s)=\frac{e^{-X_{n}}}{1+e^{-X_{n}}-s},\,n=1,2,...,  \label{geom}
\end{equation*}%
under Condition \ref{B1} restricted to $\alpha <2$ and $|\psi |<1$ have
shown that similarly to (\ref{key}) one has 
\begin{equation}
\mathbb{P}(T=n)\sim \Theta \mathbb{P}(T^{-}=n)\quad \text{as }n\rightarrow
\infty ,  \label{T1.1}
\end{equation}%
where $\Theta >0$ and 
\begin{equation*}
T^{-}=\min \{k\geq 1:S_{k}<0\}.  \label{def T^{-}}
\end{equation*}

\subsubsection{Subcritical case}

\paragraph{Weakly subcritical case.}

The best annealed results up to this moment for weakly subcritical BPRE were
obtained by Afanasyev, Boinghoff, Kersting, and Vatutin \cite{ABKV10}. It
was proved that under Conditions \ref{A2ws} and \ref{A3ws} (see below) there
exist numbers $0<\kappa ,\kappa ^{\prime }<\infty $ such that 
\begin{equation*}
\mathbb{P}\left( Z_{n}>0\right) \ \sim \ \kappa \ \mathbb{P}\left( \min
(S_{1},\ldots ,S_{n})\geq 0\right) \sim \ \kappa ^{\prime }\frac{\gamma ^{n}%
}{na_{n}}\ ,
\end{equation*}%
where $\gamma $ is from (\ref{dgamm}) and $a_{n}$ are from (\ref{defcoefal}%
). Observe that this result is similar to that in the critical case (\ref%
{ratio}), whereas it is no longer true in the intermediately and strongly
subcritical cases (see e.g.,~\cite{GKV}).

\begin{condition}
\label{A3ws} For some positive $\varepsilon$ and $a$ 
\begin{equation*}
\mathbf{E}[(\ln ^{+}\zeta (a))^{\alpha +\varepsilon }]\ <\ \infty .
\label{momcowscr}
\end{equation*}
\end{condition}

\paragraph{Intermediately subcritical case.}

Afanasyev, B\"{o}inghoff, Kersting, and Vatutin \cite{ABKV11} proved that
given $\beta =1$ and Conditions \ref{A2ws}, \ref{A3ws}  there are a constant 
$0<\theta <\infty $ and a sequence $l_{n}$ slowly varying at infinity such
that 
\begin{equation*}
\mathbb{P}\left( Z_{n}>0\right) \ \sim \ \theta \gamma ^{n}\ \mathbb{P}%
\left( S_{n}<\min (S_{1},\ldots ,S_{n-1})\right) \sim \frac{\gamma ^{n}l_{n}%
}{n^{1-1/\alpha }},\quad n\rightarrow \infty .
\end{equation*}%
Under a slightly different assumption the last relation was first proved in~%
\cite{V04}.

\paragraph{Strongly subcritical case.}

Guivarc'h and Liu in \cite{GuLi01} have shown for strongly subcritical
processes satisfying 
\begin{equation}  \label{A2ss}
\mathbb{E}\left[ \xi \ln ^{+}\xi \right] <\infty
\end{equation}
that for some $0<c\leq 1$ 
\begin{equation}
\mathbb{P}\left( Z_{n}>0\right) \thicksim c\left( \mathbb{E}\left[ \xi %
\right] \right) ^{n},\;n\rightarrow \infty.  \label{solo}
\end{equation}%
This asymptotic formula was originally established by D'Souza and Hambly in 
\cite{DH} under an extra moment assumption (also a similar statement was
proved in \cite{DH} assuming ergodicity of the environment). For the
fractional-linear case the statement (\ref{solo}) was obtained in \cite{af98}%
.

Condition \eqref{A2ss} is a counterpart of the classical $xlogx$ condition
for supercritical and subcritical branching processes in constant
environment (see \cite{AN}). It holds, in particular, under the restriction %
\eqref{bounded}.

\subsection{\protect\bigskip Conditional limit theorems}

Given that a BPRE has survived during a long period of time, what can be
said about its demographic history? The theorems of this section reveal
different successful survival strategies adjusted to the variable
reproduction strength of the population living in a random environment.

\subsubsection{\protect\bigskip Critical case}

Afanasyev, Geiger, Kersting, and Vatutin \cite{4h} proved so far the best
annealed conditional limit theorem. According to this theorem a surviving
critical BPRE behaves in a `supercritical' manner. Supercritical branching
processes (whether classical or in random environment) grow exponentially
fast $Z_{n}e^{-S_{n}}\rightarrow W\;a.s.$, where $W$ is a non-degenerate
random variable (see \cite{AN}, \cite{Sevast}).

For integers $0\leq r\leq n$ let $X_{t}^{r,n},t\in \lbrack 0,1],$ be the
rescaled generation size process given by 
\begin{equation}
X_{t}^{r,n}\ =\ Z_{r+(n-r)t}\cdot e^{-S_{r+(n-r)t}}\ ,\quad t\in \lbrack
0,1]\ .  \label{process}
\end{equation}%
In \cite{4h} it was proved that if \eqref{dpar} and Condition \ref{A2} (or
Condition \ref{B1} with $|\psi |<1$ and Condition \ref{B2}) are valid, then, 
\begin{equation}
\mathcal{L}\big(X_{t}^{r_n,n},\,t\in \lbrack 0,1]\,\big|\,Z_{n}>0\big)%
\overset{\mathcal{D}}{\rightarrow }\mathcal{L}\big(W_{t},t\in \lbrack 0,1]%
\big)\ \mbox{ as }n\rightarrow \infty ,  \label{llimit}
\end{equation}%
where the symbol $\overset{\mathcal{D}}{\rightarrow }$ stands for the
convergence in distribution in the Skorokhod topology in space $D[0,1]$ of c%
\`{a}dl\`{a}g functions on the unit interval or in other spaces which we
will meet later on. Here $r_{1},r_{2},\ldots $ is a sequence of positive
integers such that $r_{n}\leq n$ and $r_{n}\rightarrow \infty$. This
sequence is introduced to exclude the early part of the population history.

The limiting process $\left\{ W_{t},t\in \lbrack 0,1]\right\} $ is a
stochastic process with a.s. constant paths, i.e., $\mathbb{P}\left( W_{t}=W%
\text{ for all }t\in \lbrack 0,1]\right) =1$ for some positive random
variable $W$. We see that the growth of the BPRE is mainly determined by the
ARW, namely by the sequence $(e^{S_{n}})_{n\geq 0}$, with the fine structure
of the random environment being summarized by $W$. Afanasyev~\cite{Af3} \ in
the case $r_{n}=0$ established the convergence result (\ref{llimit}) for $%
t\in (0,1)$ under a stronger Condition \ref{A1(a)}.

It was also shown in \cite{4h} that given Condition \ref{B1} with $|\psi |<1$
and Condition~\ref{B2} there exists a slowly varying sequence $%
l_{1},l_{2},\ldots $ such that%
\begin{equation}
\mathcal{L}\Big(n^{-\frac{1}{\alpha }}l_{n}\,\ln Z_{nt},\,t\in \lbrack 0,1]\,%
\Big|\;Z_{n}>0\Big)\overset{\mathcal{D}}{\rightarrow }\mathcal{L}%
(L_{t}^{+},\,t\in \lbrack 0,1])\ \mbox{ as }n\rightarrow \infty ,
\label{funcT1}
\end{equation}%
where $L_{t}^{+},\,t\in \lbrack 0,1],$ denotes the meander of a strictly
stable process with index $\alpha $, which is a strictly stable L\'{e}vy
process conditioned to stay positive on the time interval $(0,1]$ (see \cite%
{Do85} and \cite{Dur78}). Convergence (\ref{funcT1}) was earlier established
by Kozlov \cite{K99} under the second moment assumption on the increments of
the ARW.

\bigskip Afanasyev \cite{Af2} studied the behavior of the process up to the
moment of extinction and proved that for a critical BPRE satisfying
Condition \ref{A1(a)} and some additional assumptions, as $n\rightarrow
\infty ,$ 
\begin{equation*}
\mathcal{L}\left( \frac{1}{\sigma \sqrt{n}}\ln (Z_{nt}+1),\,t\in \lbrack
0,\infty )\,\Big |\,Z_{n}>0\right) \overset{\mathcal{D}}{\rightarrow }%
\mathcal{L}\left( \mathbb{W}_{t}^{+},\,t\in \lbrack 0,\infty )\right) .
\label{funcT}
\end{equation*}%
The limiting process $\mathbb{W}_{t}^{+},\,t\in \lbrack 0,\infty ),$ is a
Brownian meander $W_{t}^{+}$ $\ $for $t\leq 1$ continued by a "stopped"
Brownian motion $W_{t},t\geq 1$ with $W_{1}=W_{1}^{+}$ killed at the moment 
\begin{equation*}
\tau _{0}=\inf \{t\geq 1:W_{t}=0\}.
\end{equation*}

\bigskip B\"{o}inghoff, Dyakonova, Kersting, and Vatutin \cite{BDKV}
consider a linear-fractional BPRE satisfying Conditions \ref{A1(a)}, \ref{A0}%
, and Condition \ref{B2} with $\alpha=2$. They have shown that 
\begin{equation}
\mathcal{L}\left( Z_{n}\,|\,T=n+1\right) \overset{d}{\rightarrow }\mathcal{L}%
\left( Y\right) ,\qquad n\rightarrow \infty ,  \label{DiskDist}
\end{equation}%
where $Y$ is a positive integer-valued random variable being finite with
probability 1, and furthermore, that for any $\delta \in \left( 0,1/2\right) 
$ as $n\rightarrow \infty $ 
\begin{equation*}
\mathcal{L}\left( Z_{nt}e^{-S_{nt}},t\in \left[ \delta ,1-\delta \right] \,%
\Big|\,T=n+1\right) \overset{\mathcal{D}}{\rightarrow }\mathcal{L}\left(
W_{t},t\in \left[ \delta ,1-\delta \right] \right).
\end{equation*}%
Here the limiting process is a stochastic process with a.s. constant paths,
i.e., $\mathbb{P}(W_{t}=W$ for all $t\in \lbrack 0,1])=1$ for some positive
random variable $W$, and convergence in distribution holds with respect to
the Skorokhod topology in the space $D[\delta ,1-\delta ]$ of cadlag
functions on the interval $\left[ \delta ,1-\delta \right] $.

Vatutin and Wachtel \cite{VW09} have established for the critical BPREs
whose increments of the ARWs belong to $\alpha $-stable law with $\alpha <2$
the following phenomenon of sudden extinction. Given the process dies at a
remote time $T=n+1$, the log-size $\ln Z_{nt}$ grows roughly as $n^{1/\alpha
}$ up to moment $n$ and then the process instantly dies out, so that $\ln
Z_{n}$ is of order $n^{1/\alpha }$ while $Z_{n+1}=0$. This may be
interpreted as the evolution of the process in a favorable environment up to
moment $n$ stopped by a sudden extinction of the population at moment $T=n+1$
due to a "catastrophic" change of the environment. Note that this phenomenon
is in a sharp contrast with the case of finite variance (see (\ref{DiskDist}%
)) where the distribution of the number of individuals in the population
just before the moment of its extinction converges to a distribution
concentrated on the set of positive integers.

Let us formulate the results of \cite{VW09} precisely. According to Durrett 
\cite{Dur78} under the Condition \ref{B1} there exists a nonnegative random
variable $\Lambda $ such that for $a_{n}$ from (\ref{defcoefal}) 
\begin{equation*}
\lim_{n\rightarrow \infty }\mathbb{P}(S_{n}\leq xa_{n}|S_{1}>0,\ldots
,S_{n}>0)=\mathbb{P}(\Lambda \leq x)\quad \text{for all }x\geq 0.
\label{meander}
\end{equation*}%
It is proved in \cite{VW09} that if Condition \ref{B1} holds restricted to $%
\alpha <2$ and $|\psi |<1$ together with Condition \ref{B2}, then for every $%
x>0$ 
\begin{equation*}
\lim_{n\rightarrow \infty }\frac{\mathbb{P}(Z_{n-1}>e^{xa_{n}};T=n)}{\mathbb{%
P}(T^{-}=n)}=\Theta \frac{\mathbb{E}\left[ \Lambda ^{-\alpha };\Lambda >x%
\right] }{\mathbb{E}\left[ \Lambda ^{-\alpha }\right] },  \label{T3.2}
\end{equation*}%
where (compare with \eqref{T1.1}) 
\begin{equation*}
\Theta :=\liminf_{n\rightarrow \infty }\frac{\mathbb{P}(T=n)}{\mathbb{P}%
(T^{-}=n)}.
\end{equation*}%
Moreover, it is found that in the geometric case 
\begin{equation*}
\lim_{n\rightarrow \infty }\mathbb{P}(Z_{n-1}>e^{xa_{n}}|T=n)=\frac{\mathbb{E%
}\left[ \Lambda ^{-\alpha };\Lambda >x\right] }{\mathbb{E}\left[ \Lambda
^{-\alpha }\right] }.  \label{T4}
\end{equation*}

\subsubsection{Subcritical case}

\label{S322}

\paragraph{Weakly subcritical case.}

Afanasyev, B\"{o}inghoff, Kersting, and Vatutin \cite{ABKV10} proved for
weakly subcritical processes that under Conditions \ref{A2ws} and \ref{A3ws}
the conditional laws $\mathcal{L}(Z_{n}\,|\,Z_{n}>0)$, $n\geq 1$, converge
weakly to some probability distribution on the natural numbers. Moreover, it
was shown that the sequence $\mathbb{E}[Z_{n}^{\vartheta }\,|\,Z_{n}>0]$ is
bounded for any $\vartheta <\beta $, implying convergence to the
corresponding moments of the limit distribution. It was also shown that for $%
X_{t}^{r,n}$ defined in (\ref{process}) 
\begin{equation*}
\mathcal{L}\big(X_{t}^{r_n,n},t\in \lbrack 0,1]\ \big|\ Z_{n}>0\big)\ 
\overset{\mathcal{D}}{\rightarrow }\mathcal{L}\big(W_{t},\ t\in \lbrack 0,1]%
\big),\ n\rightarrow \infty ,
\end{equation*}%
weakly in the Skorokhod space $D[0,1]$, where $r_1,r_2,\ldots$ are natural
numbers such that $r_{n}<n/2$ and $r_{n}\rightarrow \infty .$ Here $W_{t}=W$
a.s. for all $t\in \lbrack 0,1]$ and $\mathbb{P}\left( 0<W<\infty \right) \
=\ 1$. Compared to the critical case result \eqref{llimit} now both the
initial and the most recent part of the population history are excluded.
Earlier versions of this results can be found in \cite{af98} and \cite{GKV}.

Thus, we have the following pattern characterising the weakly subcritical
case: given $Z_{n}>0$ the value of $Z_{k}$ are of bounded order for $k$
close to $0$ and close to $n$. For $1\ll k\ll n$ the demographic dynamics of 
$Z_{k}$ follows the value of the quenched conditional mean $\mathcal{E}\left[
Z_{k}\right] =e^{S_{k}}$ in a completely deterministic manner up to a random
factor $W>0$, resembling the behavior of supercritical branching processes.

\paragraph{ Intermediately subcritical case.}

For intermediately subcritical BPREs under the annealed approach Afanasyev, B%
\"{o}inghoff, Kersting, and Vatutin \cite{ABKV11} proved that under
Conditions \ref{A2ws} and \ref{A3ws} the unscaled population size $Z_{n}$
conditioned on $\left\{ Z_{n}>0\right\} $ converges in distribution as $%
n\rightarrow \infty $. Under slightly different assumption this result was
first obtained by Vatutin in \cite{V04}.

Denote by $\hat{L}_t$ the strictly stable process $L_t$ with parameter $%
\alpha$ conditioned to have its minimum at time $1$. Let $e_{1},e_{2},\ldots 
$ denote the excursion intervals of $\hat{L}_t$ between consecutive local
minima and put $j\left( t\right) =i$ for $t\in e_{i}$. In \cite{ABKV11}\ it
was also proved that Conditions \ref{A2ws} and \ref{A3ws} imply for $%
0<t_{1}<t_{2}<...<t_{k}<1$, as $n\rightarrow \infty ,$%
\begin{align*}
\mathcal{L}& \left( \frac{Z_{nt_{1}}}{\exp \left( S_{nt_{1}}-\min_{k\leq
nt_{1}}S_{k}\right) },\ldots ,\frac{Z_{nt_{k}}}{\exp \left(
S_{nt_{k}}-\min_{k\leq nt_{k}}S_{k}\right) }\ \big|\ Z_{n}>0\right) \\
& \quad \quad \quad \quad \quad \quad \quad \quad \quad \quad \quad \quad
\quad \quad \quad \quad \quad \quad \quad \quad \quad \quad \overset{d}{%
\rightarrow }\mathcal{L}\left( V_{j\left( t_{1}\right) },\ldots ,V_{j\left(
t_{k}\right) }\right) ,
\end{align*}%
where $V_{1},V_{2},\ldots $ are i.i.d. copes of some strictly positive
random variable $V.$ This means that if $t_{i}$ and $t_{k}$ belong to one
and the same excursion, then $j\left( t_{i}\right) $ and $j\left(
t_{k}\right) $ coincide, i.e. $V_{j\left( t_{i}\right) }=V_{j\left(
t_{k}\right) }$ with probability $1.$ Notice that here $j(t_{1})\leq
j(t_{2})\leq \ldots \leq j(t_{n})$. For the linear fractional case this
result was earlier obtained by Afanasyev in \cite{af01sc}.

Recall that $\tau \left( nt\right) $ is the time of the first minimum of the
ARW up to time $\left[ nt\right] $. Complementing the previous results it
was shown in \cite{ABKV11} that at the times of consecutive ARW minima the
population sizes have discrete limit distributions 
\begin{equation*}
\mathcal{L}\big(\left( Z_{\tau \left( nt_{1}\right) },\ldots ,Z_{\tau \left(
nt_{k}\right) }\right) \ \big|\ Z_{n}>0\ \big)\overset{d}{\rightarrow }%
\mathcal{L}\left( Y_{j\left( t_{1}\right) },\ldots ,Y_{j\left( t_{k}\right)
}\right) ,
\end{equation*}%
where $Y_{1},Y_{2},\ldots $ are i.i.d. copies of a random variable $Y$
taking values in $\mathbf{N}.$

\paragraph{Strongly subcritical case.}

Geiger, Kersting, and Vatutin \cite{GKV} studied the strongly subcritical
BPRE under condition \eqref{A2ss} and proved that 
\begin{equation*}
\lim_{n\rightarrow \infty }\mathbb{P}\left( Z_{n}=z\ \big|\ Z_{n}>0\right)
=r_{z},\;z\in \mathbf{N},
\end{equation*}%
where 
\begin{equation*}
\sum_{z=1}^{\infty }r_{z}=1\text{ and }m_{r}:=\sum_{z=1}^{\infty
}zr_{z}<\infty .
\end{equation*}

Afanasyev, Geiger, Kersting, and Vatutin \cite{AGKV2} have introduced the
following condition in the strongly subcritical case.

\begin{condition}
\label{A2'ss} Suppose that (see \eqref{ddefeta}) 
\begin{equation*}
\mathbb{E}\left[ e^{X}\ln ^{+}\eta \right] <\infty .
\end{equation*}
\end{condition}

Since in the strongly subcritical case $\mathbb{E}\left[ e^{X}\right]
<\infty $, Condition \ref{A2'ss} holds in particular if the random offspring
distribution \eqref{Q} has uniformly bounded support. It also holds if %
\eqref{Q} is a Poisson distribution with random mean, so that $\eta =1$
a.s., or if \eqref{Q} is a geometric distribution on $\mathbf{N}_{0} $ where 
$\eta =2$ \ a.s.

It was shown in \cite{AGKV2} that with 
\begin{equation*}
0=:i_{n,0}<i_{n,1}<i_{n,2}<...<i_{n,k}<i_{n,k+1}:=n
\end{equation*}%
the following weak convergence holds 
\begin{equation*}
\mathcal{L}\big(\left( Z_{j}\right) _{0\leq j\leq m},\left(
Z_{i_{n,1}+j}\right) _{0\leq j\leq m},...,\left( Z_{i_{n,k}+j}\right)
_{0\leq j\leq m},\left( Z_{n-j}\right) _{0\leq j\leq m}\ \big|\ Z_{n}>0\big)%
\ 
\end{equation*}%
\begin{equation*}
\overset{d}{\rightarrow }\mathcal{L}_{\delta _{1}}\ \big(\left( Y_{j}\right)
_{0\leq j\leq m}\ \big)\otimes \mathcal{L}_{\widehat{r}}\ \big(\left(
Y_{j}\right) _{0\leq j\leq m}\ \big)^{\otimes k}\otimes \mathcal{L}_{r}\ %
\big((\widetilde{Y}_{j})_{0\leq j\leq m}\ \big),
\end{equation*}%
for every $k,m\in \mathbf{N}_{0}$ as $\min_{0\leq l\leq k}\left(
i_{n,l+1}-i_{n,l}\right) \rightarrow \infty $\ and $n\rightarrow \infty $.
Here, $\mathcal{L}_{\nu }\ \big(\left( Y_{j}\right) _{j\geq 0}\big)$ denotes
the law of the Markov chain $\left( Y_{j}\right) _{j\geq 0}$ with initial
distribution $\nu $ and transitional probabilities%
\begin{equation*}
\widehat{P}_{yz}=\frac{z\mathbb{P}\left( \xi _{1}+\ldots +\xi _{y}=z\right) 
}{y\mathbb{E}\left[ \xi \right] },\quad y,z\in \mathbf{N},
\end{equation*}%
where $\xi _{1},\ldots ,\xi _{y}$ are independent copies of the offspring
number $\xi $. The stationary distribution of $\left( Y_{j}\right) _{j\geq
0} $\ is 
\begin{equation*}
\widehat{r}_{z}=\frac{zr_{z}}{m_{r}},\;z\in \mathbf{N.}
\end{equation*}%
The process $(\widetilde{Y}_{j})_{j\geq 0}$ is the time-reversed Markov
chain with transition probabilities%
\begin{equation*}
\widetilde{P}_{yz}=\frac{\widehat{r}_{z}\widehat{P}_{zy}}{\widehat{r}_{y}}.
\end{equation*}

\section{Multi-type BPREs}

\label{Smu}

The model of multi-type BPRE was first considered by Athreya and Karlin \cite%
{AK71} and subsequently investigated by Weissener \cite{Weis}, Kaplan \cite%
{Kap}, Tanny \cite{Ta} ( see also \cite{VZ} for other references).\ Next we
describe the BPREs with $p$ types of particles $\mathbf{Z}%
_{n}=(Z_{n,1},...,Z_{n,p})$ and time $n=0,1,\ldots $ using standard notation
for the $p$-dimensional vectors:

\begin{itemize}
\item the unit vector $\mathbf{e}_{j}$ has $j$-th component equals $1$ and
all others equal zero, $j=1,\ldots ,p$,

\item all zero and all one vectors $\mathbf{0}=(0,\ldots,0)$, $\mathbf{1}%
=(1,\ldots,1)$,

\item for $\mathbf{x}=(x_{1},...,x_{p})$ and $\mathbf{y}=(y_{1},...,y_{p})$
set%
\begin{equation*}
\left\Vert \mathbf{x}\right\Vert =\sum_{i=1}^{p}|x_{i}|,\quad (\mathbf{x},%
\mathbf{y})=\sum_{i=1}^{p}x_{i}y_{i},\quad \mathbf{x}^{\mathbf{y}%
}=\prod\limits_{i=1}^{p}x_{i}^{y_{i}}.
\end{equation*}
\end{itemize}

A particle of type $i$ in a multi-type branching process can produce $z_{1}$
particles of type 1, $z_{2}$ particles of type 2, $\ldots $, and $z_{p}$
particles of type $p$ with a probability $Q_{\mathbf{z}}^{(i)}$, $\mathbf{z}%
=(z_{1},\ldots ,z_{p})$. Under the random environment assumption the
probabilities $Q_{\mathbf{z}}^{(i)}$ are random variables making random the
corresponding PGFs 
\begin{equation}
F^{(i)}(s_{1},\ldots ,s_{p})=\sum_{z_{1}=0}^{\infty }\ldots
\sum_{z_{p}=0}^{\infty }Q_{\mathbf{z}}^{(i)}\mathbf{s}^{\mathbf{z}},\
i=1,\ldots ,p.  \label{Defff}
\end{equation}%
The random PGFs for the consecutive generations $%
(F_{n}^{(1)},...,F_{n}^{(p)})_{n\geq 1}$ representing the environmental
history determine the conditional mean offspring numbers as random $p\times
p $ matrices $\mathbf{M}_{n}$ with elements $\partial F_{n}^{(i)}(\mathbf{1}%
)/\partial s_{j}$. These are IID matrices having the same distribution as 
\begin{equation}
\mathbf{M}=(M_{ij})_{i,j=1}^{p},\ M_{ij}=\frac{\partial F^{(i)}(\mathbf{1})}{%
\partial s_{j}}.  \label{matsr1}
\end{equation}

We will assume that elements $M_{ij}$ are all positive meaning that each of $%
p$ types of particles can produce any other type of particles in the next
generation. We denote by $R$ and $R_{n}$ the Perron roots, i.e. the maximal
(in absolute value) eigenvalues, for the matrices $\mathbf{M}$ and $\mathbf{M%
}_{n}$. There exist left and right eigenvectors corresponding to the Perron
eigenvalue 
\begin{align*}
& \mathbf{V}=(V_{1},...,V_{p}),\ \mathbf{V}\mathbf{M}=R\mathbf{V}, \\
& \mathbf{U}=(U_{1},...,U_{p}),\ \mathbf{M}\mathbf{U}^{\prime }=R\mathbf{U}%
^{\prime },
\end{align*}%
both with positive components and unique up to scaling. Although vectors $%
\mathbf{V}$ and $\mathbf{U}$ are random in general, throughout this section
we assume that one of the following basic conditions is valid.

\begin{condition}
\label{A0lm} There is a non-random strictly positive vector $\mathbf{v}$
such that $\left\Vert \mathbf{v}\right\Vert =1$ and 
\begin{equation*}
\mathbb{P}(\mathbf{V}=\mathbf{v})=1.
\end{equation*}
\end{condition}

\begin{condition}
\label{A0rm} There is a non-random strictly positive vector $\mathbf{u}$
such that $\left\Vert \mathbf{u}\right\Vert =1$ and 
\begin{equation*}
\mathbb{P}(\mathbf{U}=\mathbf{u})=1.
\end{equation*}
\end{condition}

The important gain of this restriction is that we get either 
\begin{equation*}
\mathbf{v}\mathbf{M}_{1}\ldots \mathbf{M}_{n}=R_{1}\cdots R_{n}\mathbf{v}
\end{equation*}%
or 
\begin{equation*}
\mathbf{M}_{1}\ldots \mathbf{M}_{n}\mathbf{u}^{\prime }=R_{1}\cdots R_{n}%
\mathbf{u}^{\prime }.
\end{equation*}%
In both cases putting $X:=\ln R,$ $X_{n}:=\ln R_{n},n\geq 1,$ we can again
introduce an ARW by $S_{n}=X_{1}+...+X_{n}$. We will refer to the conditions
on $S_{n}$ stated for the single type case even in the multi-type setting
(under Conditions \ref{A0lm} or \ref{A0rm}). In particular, the single type
classification for BPREs extends straightforwardly to this multi-type case.
For example under Condition \ref{A0lm} as $n\rightarrow \infty $ 
\begin{align*}
\sum_{i=1}^{p}v_{i}\mathcal{P}(\mathbf{Z}_{n}\neq \mathbf{0\,}|\,\mathbf{Z}%
_{0}=\mathbf{e}_{i})& \leq \min_{1\leq k\leq n}|\mathbf{vM}_{1}\cdot \cdot
\cdot \mathbf{M}_{k}|=\min_{1\leq k\leq n}R_{1}\cdots R_{k} \\
& =\exp \{\min_{1\leq k\leq n}S_{k}\}\rightarrow 0
\end{align*}%
almost surely. This implies that in the critical and subcritical cases 
\begin{equation*}
\mathcal{P}(\mathbf{Z}_{n}\neq \mathbf{0\,}|\,\mathbf{Z}_{0}=\mathbf{e}%
_{i})\rightarrow 0
\end{equation*}%
almost surely for any $i=1,\ldots ,p$. The same conclusion follows from
Condition~\ref{A0rm} as well.

\subsection{\protect\bigskip The survival probability of the multi-type BPREs%
}

Dyakonova and Vatutin have managed to extend some of the single type
annealed results to the multi-type BPREs in the critical and subcritical
cases. Some of the results mentioned in this section bring a stronger
version of the known single type results.

\subsubsection{\protect\bigskip Critical multi-type case}

Dyakonova \cite{Dyak} considered a multi-type critical BPRE satisfying %
\eqref{dpar}, Condition~\ref{A0lm}, and Condition~\ref{A2}, where $\zeta (a)$
is now defined as 
\begin{equation}
\zeta (a):=R^{-2}\sum_{\mathbf{z}\in \mathbf{N}_{a}^{p}}\sum_{i=1}^{p}v_{i}%
\sum_{j,k=1}^{p}Q_{\mathbf{z}}^{(i)}z_{j}z_{k},  \label{za}
\end{equation}%
with $\mathbf{N}_{a}^{p}$ being the set of vectors $\mathbf{z}$ with
non-negative integer-valued components such that at least one component is
larger or equal $a$. Here $R$ is the (random) Perron root of the random
matrix $\mathbf{M}$ related to $Q_{\mathbf{z}}^{(i)}$ through \eqref{Defff}
and \eqref{matsr1}. Additionally the following restriction (due to \cite%
{SpitkKe}) is imposed.

\begin{condition}
\label{A0dopclas} There exists a number $0<d<1$ such that 
\begin{equation*}
d\leq \frac{M_{i_{1},j_{1}}}{M_{i_{2},j_{2}}}\leq d^{-1},1\leq
i_{1},i_{2},j_{1},j_{2}\leq p.
\end{equation*}
\end{condition}

Under these conditions it was proved in \cite{Dyak} that (compare with %
\eqref{ratio}) 
\begin{equation}
\mathbb{P}(\mathbf{Z}_{n}\neq \mathbf{0}\,|\,\mathbf{Z}_{0}=\mathbf{e}%
_{i})\sim c_{i}n^{-(1-\rho )}l_{n},\text{ \ }i=1,\ldots ,p,  \label{qqqq}
\end{equation}%
where $c_{i}$ are positive constants, $\rho $ is the constant from the
Spitzer condition applied to $S_{n}=\ln R_{1}+...+\ln R_{n}$, and $l_{n}$ is
a sequence slowly varying at infinity. Earlier the same asymptotics (\ref%
{qqqq}) was established in \cite{Dyak9} in the linear fractional case under
Condition \ref{A1(a)} with $|\psi |<1$ and Condition \ref{A0rm}.

In \cite{Dyak} \ it was proved also a quenched type result claiming the
almost sure convergence 
\begin{equation}
\frac{\mathcal{P}(\mathbf{Z}_{n}\neq \mathbf{0\,}|\,\mathbf{Z}_{0}=\mathbf{e}%
_{i})}{\sum_{j=1}^{p}v_{j}\mathcal{P}(\mathbf{Z}_{n}\neq \mathbf{0\,}|\,%
\mathbf{Z}_{0}=\mathbf{e}_{j})}\rightarrow \mathbf{U}^{\ast },\qquad
i=1,\ldots ,p,  \label{fTh2mc}
\end{equation}%
where $\mathbf{U}^{\ast }=(U_{1}^{\ast },...,U_{p}^{\ast })$ is a random
vector such that 
\begin{equation*}
(\mathbf{v},\mathbf{U}^{\ast })=1,\;d\leq U_{i}^{\ast }\leq 1/\min
(v_{1},...,v_{p}),\qquad i=1,\ldots ,p.  \label{con1}
\end{equation*}
Furthermore, it was shown in \cite{Dyak} that the above mentioned results
remain true for an alternative set of conditions involving the counterparts
of the single type Conditions \ref{B1} and \ref{B2}.

Vatutin and Dyakonova \cite{VD2010} studied a multi-type critical BPRE whose
ARW satisfies the Spitzer-Doney condition and the mean matrices of the
reproduction laws have a common positive right eigenvector (see Condition %
\ref{A0rm}). Also in \cite{VD2010} it is assumed that the following
assumption holds.

\begin{condition}
\label{A2rm} By this condition we exclude a possibility for the defining
PGFs to take the linear form 
\begin{equation*}
F^{(i)}(\mathbf{s})=Q_{\mathbf{0}}^{(i)}+Q_{\mathbf{e}_{1}}^{(i)}s_{1}+%
\ldots +Q_{\mathbf{e}_{p}}^{(i)}s_{p}.
\end{equation*}
\end{condition}

Another condition required in \cite{VD2010} concerns the moments for the
vectors $(\xi _{i1},\ldots ,\xi _{ip})$ of the offspring numbers having
quenched distributions $\{Q_{\mathbf{z}}^{(i)}\}$ and quenched means $%
(M_{i1},\ldots ,M_{ip})$. For a given $\beta >0$ set 
\begin{equation*}
\Delta _{ij}(\beta )=\mathcal{E}\left\vert \xi _{ij}-M_{ij}\right\vert
^{\beta },\quad \Delta _{\beta }=\max_{i,j}\Delta _{ij}(\beta )
\end{equation*}%
and denote $\zeta _{\beta }:=e^{-\beta X}\Delta _{\beta }$.

\begin{condition}
\label{A3rm} There exist $\beta \in (1,2]$ and $\varepsilon >0$ such that 
\begin{equation}
\mathbb{P}(\Delta _{\beta }<\infty )=1,  \label{Asinfin}
\end{equation}%
and%
\begin{equation*}
\mathbb{E}[\ln ^{+}\zeta _{\beta }]^{1/\rho +\varepsilon }<\infty ,\quad 
\mathbb{E}[V(X)(\ln ^{+}\zeta _{\beta })^{1+\varepsilon }]<\infty .
\label{Ass22}
\end{equation*}
\end{condition}

In \cite{VD2010} the asymptotic result \eqref{qqqq} was established under %
\eqref{dpar} and Conditions \ref{A0rm}, \ref{A2rm}, and \ref{A3rm}. The
results of \cite{VD2010} generalize and extend not only the above mentioned
results in \cite{Dyak} but also the corresponding single type statements
from \cite{4h}. Indeed, in view of condition (\ref{Asinfin}) the asymptotic
relation (\ref{qqqq}) may be valid for $p=1$ if, for instance,%
\begin{equation*}
\mathbb{P}\left(\mathcal{E}\left[ \xi ^{2}1_{\{\xi \geq a\}}\right] =\infty
\right) =1  \label{21}
\end{equation*}%
for any $a\geq 0$ making in this case the conditions of the last result
weaker than those required in \cite{4h} for \eqref{ratio}.

As in the single type case, denote by $T$ the extinction moment for the
process $\mathbf{Z}_n$. It was shown by Dyakonova \cite{Dy02} that in the
multi-type linear fractional case under Conditions \ref{A1(a)}, \ref{A0lm}, %
\ref{A0rm}, \ref{A0dopclas} and some extra restrictions there exists a
constant $0<c<\infty $ such that (compare with \eqref{VD97(1.7)}) for any $%
i=1,\ldots ,p$ 
\begin{equation*}
\mathbb{P}(T=n|\mathbf{Z}_{0}=\mathbf{e}_{i})\sim cu_{i}n^{-3/2},\qquad
n\rightarrow \infty ,  \label{5}
\end{equation*}%
where $\mathbf{u}=(u_{1},...,u_{p})$ is defined in Condition \ref{A0rm}.

\subsubsection{Subcritical multi-type case}

Dyakonova \cite{Dy08} investigated the asymptotics of the survival
probability of multi-type subcritical BPREs. She has proved, that if
Conditions \ref{A0lm} and \ref{A0dopclas} are valid and $\mathbb{E}\left[ X%
\right] <0$, then \eqref{fTh2mc} holds $\mathbb{P}$-a.s., as $n\rightarrow
\infty $. Furthermore, in \cite{Dy08} it was shown that in the strongly
subcritical case as $n\rightarrow \infty $ (compare with \eqref{solo}) 
\begin{equation*}
\mathbb{P}(\mathbf{Z}_{n}\neq \mathbf{0}\,|\,\mathbf{Z}_{0}=\mathbf{e}%
_{i})\sim c_{i}(\mathbb{E}\left[ R\right] )^{n},\text{ \ }\
c_{i}>0,\;i=1,\ldots ,p  \label{sttth}
\end{equation*}%
under Conditions \ref{A0lm}, \ref{A0dopclas}, and \ref{A2'ss}, with $\eta
=\zeta (0)$ defined by \eqref{za}.

\subsection{ Functional limit theorem in the critical case}

Vatutin and Dyakonova \cite{VD2010} have assumed that assumption \eqref{dpar}
is valid together with Conditions \ref{A0rm}, \ref{A2rm}, and Condition \ref%
{A3rm} with $\beta =2$. They considered a family of the processes 
\begin{equation*}
W_{r,n}(t):=e^{-S_{r+(n-r)t}}(\mathbf{Z}_{r+(n-r)t},\mathbf{u}),\text{ }t\in
\lbrack 0,1]
\end{equation*}%
for any given pair of integers $0\leq r\leq n$. Letting $r_{1},r_{2},\ldots $
be a sequence of positive integers such that $r_{n}\leq n$ and $%
r_{n}\rightarrow \infty $, it was proved that as $n\rightarrow \infty ,$%
\begin{equation}
\mathcal{L}\left( W_{r_{n},n}(t),t\in \lbrack 0,1]\,|\,\mathbf{Z}_{n}\neq 
\mathbf{0;Z}_{0}=\mathbf{z}\right) \overset{\mathcal{D}}{\rightarrow }%
\mathcal{L}\left( W_{\mathbf{z}}(t),t\in \lbrack 0,1]\right) ,
\label{fTh1.3m}
\end{equation}%
where the limiting process is a stochastic process with a.s. constant paths,
i.e., $\mathbb{P}(W_{\mathbf{z}}(t)=W_{\mathbf{z}}$ for all $t\in \left[ 0,1%
\right] )=1$ with $\mathbb{P(}0<W_{\mathbf{z}}\text{ }<\infty )=1$.
Convergence (\ref{fTh1.3m}) is a generalization of \eqref{llimit} in two
directions. First, it is proved for the multitype case, and second, even for
the single-type case the conditions under which (\ref{fTh1.3m}) is
established in \cite{VD2010} are weaker than those used in \cite{4h}.

\section{Quenched results for the critical BPREs}

\label{Squ}

Vatutin and Dyakonova \cite{VD} and \cite{VD2} applied the quenched approach
to the critical BPREs meeting the Doney-Spitzer condition (\ref{dpar}).
Their results generalize those in \cite{VDfr} and \cite{VDWie} in which much
stronger conditions on the characteristics of the branching processes are
imposed.


\subsection{Survival probability conditioned on the environment}

The following condition involving the renewal functions (\ref{DevU}) is dual
to Condition \ref{A2} dealing with the renewal function (\ref{DefV}).

\begin{condition}
\label{A2vd} There are numbers $\varepsilon >0$ and $a\in \mathbf{N}_{0}$
such that 
\begin{equation*}
\mathbb{E}[(\ln ^{+}\zeta (a))^{\frac{1}{1-\rho }+\varepsilon }]<\infty
\quad \mbox{and}\quad \mathbb{E}[U(-X)(\ln ^{+}\zeta (a))^{1+\varepsilon
}]<\infty ,  \label{A22}
\end{equation*}%
where $\zeta (a)$ is from (\ref{defdzet}).
\end{condition}

It was proved in \cite{VD2} that under \eqref{dpar} and Conditions \ref{A2}, %
\ref{A2vd} 
\begin{equation}
e^{-S_{\tau \left( n\right) }}\mathcal{P}\left( Z_{n}>0\right)\overset{d}{%
\rightarrow }\zeta ,\qquad n\rightarrow \infty ,  \label{formoy}
\end{equation}%
where the random variable $\zeta \in \lbrack 0,1]$ is positive with
probability~1. According to~\eqref{formoy} the asymptotic behavior of the
survival probability is mainly determined by the minimal value of the ARW
over the time interval $[0,n]$.

\subsection{Convergence of one-dimensional distributions}

Denote%
\begin{equation*}
\widehat{Y}_{n}:=\frac{Z_{n}}{\mathcal{E}\left[ Z_{n}\,|\,Z_{n}>0\right] }%
\quad  \label{mstar1}
\end{equation*}%
and 
\begin{equation*}
\mathcal{M}_{\;}^{(n)}(dx):=\mathcal{P}\left( \widehat{Y}_{n}\in
dx\left\vert Z_n>0\right. \right) .
\end{equation*}
Theorem 1 and Lemma 7 from \cite{VD2} yield that under \eqref{dpar} and
Conditions \ref{A2}, \ref{A2vd} there exists a random measure $\mathcal{M}$
(being proper and nondegenerate with probability~1) such that, as $%
n\rightarrow \infty ,$ 
\begin{align}
\int_{0}^{\infty }e^{-\lambda x}\mathcal{M}_{\;}^{(n)}(dx)& \overset{d}{%
\rightarrow }\Phi \left( \lambda \right) :=\int_{0}^{\infty }e^{-\lambda x}%
\mathcal{M}(dx),  \label{ddistcont2} \\
\int_{0}^{\infty }e^{-\lambda x}x\mathcal{M}_{\;}^{(n)}(dx)& \overset{d}{%
\rightarrow }\Psi \left( \lambda \right) :=\int_{0}^{\infty }e^{-\lambda x}x%
\mathcal{M}(dx)=-\frac{d}{d\lambda }\Phi \left( \lambda \right) .
\label{addit}
\end{align}

In the linear-fractional case the limiting Laplace transform is
deterministic $\Phi \left( \lambda \right) =(1+\lambda)^{-1}$ and represents
an exponential distribution with parameter 1. Notice that in this case $\Psi
\left( \lambda \right) =(1+\lambda )^{-2}$ corresponds to the sum of two
independent Exp$(1)$ random variables. Thus in the linear-fractional case we
get an analogue of the corresponding Yaglom-type limit theorem for the
ordinary critical Galton-Watson processes.

For the case when the distribution of $X$ is absolutely continuous
convergences (\ref{formoy}) and (\ref{ddistcont2}) were proved in~\cite%
{VDWie}.

\subsection{Convergence of finite-dimensional distributions}

For an integer $b\geq 2$ and tuples $\bar{t}=(t_{0},t_{1},\ldots ,t_{b})$, $%
\overline{\lambda }=(\lambda _{1},\ldots ,\lambda _{b})$, where $%
0=t_{0}<t_{1}<\cdots <t_{b}=1$ and $\lambda _{i}\geq 0$, $i=1,\ldots ,b$,
set 
\begin{equation*}
\Phi ^{(n)}(\bar{t},\overline{\lambda }):=\mathcal{E}\left[ \exp \left\{
-\sum_{i=1}^{b}\lambda _{i}\widehat{Y}_{nt_{i}}\right\} \,\Big|\,Z_{n}>0%
\right]  \label{matoz}
\end{equation*}%
and for a positive integer $d\leq b$ and a vector $\overline{r}%
=(r_{0},r_{1},r_{2},\ldots ,r_{d})$ with integer-valued coordinates $%
0=r_{0}<r_{1}<r_{2}<\cdots <r_{d}=b$ introduce the event 
\begin{align*}
\mathcal{U}(\bar{t};\overline{r},n):= \Big\{\tau (nt_{1})&=\cdots =\tau
(nt_{r_{1}})<\tau (nt_{r_{1}+1})=\cdots  \notag \\
&=\tau (nt_{r_{2}})<\cdots <\tau (nt_{r_{d-1}+1})=\cdots =\tau (nt_{r_{d}})%
\Big\}.  \label{deU}
\end{align*}%
It was proved in \cite{VD2} that if $X$ meets Conditions \ref{B1}, \ref{A2}, %
\ref{A2vd}, then, as $n\rightarrow \infty, $ 
\begin{equation*}
\left\{ \Phi ^{(n)}(\bar{t},\overline{\lambda })\,\big|\,\mathcal{U}\left( 
\bar{t};\overline{r},n\right) \right\} \overset{d}{\to }\prod_{i=1}^{d-1}%
\Psi_{i}\Bigg(\sum_{j=r_{i-1}+1}^{r_{i}}\lambda _{j}\Bigg)\Phi_{d}\Bigg(%
\sum_{j=r_{d-1}+1}^{b}\lambda _{j}\Bigg),
\end{equation*}%
where $\Psi_{i}(\lambda )$, $i=1,\ldots ,d-1$ and $\Phi_d(\lambda )$ are
independent random functions distributed as $\Psi (\lambda )$ in (\ref{addit}%
) and $\Phi (\lambda )$ in (\ref{ddistcont2}) respectively.

Loosely speaking, this result shows that the trajectory $(Z_0,\ldots,Z_n)$
is partitioned into independent pieces generated by strict descending ladder
moments of the ARW. The population size at the moment $j$ is proportional to 
$e^{S_{j}-S_{\tau (j)}}$ while the maximal size of the population on the
interval $[0,n]$ is proportional to $e^{\max_{0\leq j\leq n}(S_{j}-S_{\tau
(j)})}$ (with a random positive factor $\le1$).

Let $L_{t}$ be an $\alpha$-stable L\'{e}vy process. For fixed $%
0=t_{0}<t_{1}<\cdots <t_{b}=1$ put $\omega _{p}:=\inf_{t_{0}\leq u\leq
t_{p}}L_{u}$, $0\leq p\leq b$. Let $D$~be the random number of different
elements in the tuple $\omega _{0},\omega _{1},\ldots ,\omega _{b}$, and let 
$\bar{R}$ be the random $D$-dimensional vector with components 
\begin{equation*}
R_{0}=0,R_{i+1}=\max \big\{k\geq R_{i}+1:\omega _{R_{i}+1}=\omega _{k}\big\}
\end{equation*}%
so that $R_D=b$. In \cite{VD2} it was established that if $X$ satisfies
Conditions \ref{B1}, \ref{A2}, \ref{A2vd}, then as $n\rightarrow \infty $ 
\begin{equation*}
\Phi ^{(n)}(\bar{t},\overline{\lambda })\overset{d}{\to }\prod_{i=1}^{D-1}%
\Psi_{i}\Bigg(\sum_{j=R_{i-1}+1}^{R_{i}}\lambda _{j}\Bigg) \Phi_{D}\Bigg(%
\sum_{j=R_{D-1}+1}^{b}\lambda _{j}\Bigg),
\end{equation*}%
where $(\Psi_{1}(\lambda ),\Phi_1 (\lambda ),\Psi_{2}(\lambda ),\Phi_2
(\lambda ),\ldots)$ are as before and independent of the L\'{e}vy process $%
L_t$ defining the vector $\bar{R}$.

Relations (\ref{formoy}) and (\ref{ddistcont2}) imply that if $t\in (0,1]$
is \textit{fixed} then, given $Z_n>0,$ the distribution of the random
variable $Z_{nt}e^{S_{\tau (nt)}-S_{nt}}$ converges in the specified sense
to a proper distribution with no atom at zero. This means, roughly speaking,
that if the process survives up to moment $n$, then an earlier population
size $Z_{nt}$ is proportional to $e^{S_{nt}-S_{\tau (nt)}}$. Thus, contrary
to the conditional limit theorems for the classical critical or
supercritical branching processes (in which the scaling functions for the
population size increase with time either linearly or exponentially), the
scaling function for the population size in the critical BPRE is subject to
large random oscillations.

This indicates that the corresponding population passes through a number of
bottlenecks at the moments around the consecutive points of minima of the
ARW. Vatutin and Dyakonova in \cite{VDfr}, \cite{VD}, \cite{VD2}, \cite%
{VDWie}, and \cite{VD3} have investigated this phenomena in detail under the
quenched approach. They shown that the distribution of the number of
individuals in the process at the moments close to $\tau (nt),$ $t\in (0,1]$
conditioned on survival up to time $n$ converges to a discrete distribution.
Thus, in contrast to the \textit{fixed} points of the form $nt$ where the
size of the population is exponentially large (see \cite{VD2}), the size of
the population at the (random) point of global minimum of the ARW becomes
drastically small but then it grows again exponentially. This reminds the
typical demographics of real biological populations which during their
evolution have "favorable periods" (rapid growth of the population size) and
"unfavorable periods" (quick extinction when only a few representatives of
the population survive) followed by another period of rapid growth. Note a
similarity between this picture and the intermediately subcritical case
under the annealed approach (see Section~\ref{S322}).

\subsection{Discrete limit distributions}

Consider the population size near the time of the global minimum of the ARW.
Under \eqref{dpar} and Conditions \ref{A2}, \ref{A2vd} it was proved in \cite%
{VD3} that for any $m\in \mathbb{Z}$ and $t\in (0,1)$ 
\begin{align*}
\mathcal{E}\left[ s^{Z_{\tau (n)+m}}\Big|\;Z_n>0\right] & \overset{d}{%
\rightarrow }\varphi _{m}(s), \\
\mathcal{L}\left( \mathcal{E}\left[ s^{Z_{\tau (nt)+m}}\Big|\;Z_n>0\right] %
\Big|\;\tau (n)>nt\right) & \overset{d}{\rightarrow }\mathcal{L}\left(
\varphi _{m}^{\ast }(s)\right) ,
\end{align*}%
where $\varphi _{m}(s)$ and $\varphi _{m}^{\ast }(s)$ can be obtained
precisely. Observe, that the second limiting PGF does not depend on a
particular value of $t\in (0,1)$.

Next we give another two convergence results from \cite{VD3}. The first of
them characterizes the conditional distribution of the number of individuals
in the population at two sequential moments located in a vicinity of the
moment when the global minimum of the associated random walk is attained,
while the second convergence describes the same distribution in a vicinity
of the moment $\tau (nt)$ given $\tau (n)>nt.$ So, in \cite{VD3} proved that
under \eqref{dpar} and Conditions \ref{A2}, \ref{A2vd} for any $m\in \mathbb{%
Z}$ and any $t\in \left( 0,1\right) ,s_{1},s_{2}\in \lbrack 0,1),$ as $%
n\rightarrow \infty $, 
\begin{align*}
\mathcal{E}\left[ s_{1}^{Z(\tau (n)+m)}s_{2}^{Z(\tau (n)+m+1)}\Big|\;Z(n)>0%
\right] & \overset{d}{\rightarrow }\varphi _{m}(s_{1},s_{2}), \\
\mathcal{L}\left( \mathcal{E}\left[ s_{1}^{Z(\tau (nt)+m)}s_{2}^{Z(\tau
(nt)+m+1)}\Big|\;Z(n)>0\right] \Big|\;\tau (n)>nt\right) & \overset{d}{%
\rightarrow }\mathcal{L}\left( \varphi _{m}^{\ast }(s_{1},s_{2})\right) ,
\end{align*}%
where $\varphi _{m}(s_{1},s_{2})$ and $\varphi _{m}^{\ast }(s_{1},s_{2})$
can be obtained precisely.

\subsection{Conditioning on the precise time of extinction}

Critical BPREs conditioned on extinction at a given moment were investigated
by Vatutin and Kyprianou \cite{VK08} using the quenched approach. Let 
\begin{equation*}
\tau (l,n):=\min \left\{ k\in \lbrack l,n]:S_{k}=\min_{l\leq p\leq
n}S_{p}\right\}
\end{equation*}%
be the left-most point on $[l,n]$ at which the minimal value of the ARW on
the interval $[l,n]$ is attained. Assuming $\left\{ T=n\right\} $ the
authors of \cite{VK08} proved conditional limit theorems describing the
asymptotic behavior, as $n\rightarrow \infty ,$ of the distribution of the
population sizes at moments $nt,$ $t\in (0,1)$ and at moments close to $\tau
(nt)$. It turned out that if $X $ belongs to the domain of attraction of a
stable law with parameter $\alpha \in (0,2]$, then (contrary to the annealed
approach \cite{VW09}) under the quenched approach the phenomenon of sudden
extinction does not occur.

Let 
\begin{equation*}
A_{m,n}:=\frac{\mathcal{P}\left( Z_{n}>0|Z_{0}=1\right) }{\mathcal{P}\left(
Z_{n}>0|Z_{m}=1\right) }\mathcal{P}\left( Z_{n}=0|Z_{m}=1\right) b_{m},
\label{DefOO}
\end{equation*}%
where $b_{m}:=\sum_{j=0}^{m-1}\eta _{j+1}e^{-S_{j}}/2$, and $\eta _{j}$ are
from (\ref{ddefeta}). It is shown in \cite{VK08} that in the
linear-fractional case given assumption \eqref{dpar} and Conditions \ref{A2}%
, \ref{A2vd} are valid the following is true. For any $m\in \mathbb{Z}$, $%
t\in (0,1),$ $s\in (0,1]$, $\lambda \in (0,\infty )$, as $n\rightarrow
\infty $ 
\begin{equation}
\mathcal{L}\left( \mathcal{E}\left[ s^{Z_{\tau (nt)+m}}|T=n\right]
\left\vert \tau (n)\geq nt\right. \right) \overset{d}{\rightarrow }\mathcal{L%
}\left( s\left( \frac{1-\Theta _{m}}{1-\Theta _{m}s}\right) ^{2}\right) ,
\label{fTclose1}
\end{equation}%
\begin{equation}
\mathcal{L}\left( \mathcal{E}\left[ s^{Z_{\tau (nt,n)+m}}|T=n\right]
\left\vert \tau (n)<nt\right. \right) \overset{d}{\rightarrow }\mathcal{L}%
\left( s\left( \frac{1-\theta _{m}}{1-\theta _{m}s}\right) ^{2}\right) ,
\label{fTclose2}
\end{equation}%
where $\Theta _{m}\in (0,1)$ , $\theta _{m}\in (0,1)$ with probability 1,
and furthermore 
\begin{eqnarray}
\mathcal{E}\left[ \exp \left\{ -\lambda \frac{Z_{nt}}{A_{nt,n}}\right\} |T=n%
\right]\overset{d}{\rightarrow }&\frac{1}{\left( 1+\lambda \right) ^{2}}.
\label{fTfar1}
\end{eqnarray}

According to \cite{VK08} on the set $\tau (n)\geq nt$ the random variable $%
A_{nt,n}e^{S_{\tau (nt)}-S_{nt}}$ converges in distribution, as $%
n\rightarrow \infty $, to a random variable being finite and positive with
probability 1. Thus, for such moments $nt$ the normalization in (\ref{fTfar1}%
) is essentially specified by the past behavior of the ARW. On the other
hand, given $\tau (n)<nt$ the random variable $A_{nt,n}e^{S_{\tau
(nt,n)}-S_{nt}}$ converges in distribution, as $n\rightarrow \infty $, to a
random variable being finite and positive with probability 1. Thus, for such
moments $nt$ the scaling in (\ref{fTfar1}) is essentially specified by the
future behavior of the ARW.

This fact allows us to give the following non-rigorous interpretation of the
mentioned results. If the process dies out at a distant moment $T=n$, then
it happens not as a unique catastrophic event. Before the extinction moment
the evolution of the process consists of a number of "bad" periods
characterized by small population sizes. According to (\ref{fTclose1}) and (%
\ref{fTclose2}), such periods are located in the vicinities of \textit{random%
} points $\tau (nt)1_{\left\{ \tau (n)\geq nt\right\}} $ and $\tau
(nt,n)1_{\left\{ \tau (n)<nt\right\}} .$ On the other hand, at \textit{%
nonrandom }points $nt,\,t\in (0,1),$ the size of the population by (\ref%
{fTfar1}) is big. Hence $\ln Z_{nt}$ grows like $S_{nt}-S_{\tau (nt)}$ if $%
\tau (n)>nt$ and like $S_{nt}-S_{\tau (nt,n)}$ if $\tau (n)<nt$. Thus, the
process $\{Z_{tT},0\leq t\leq 1\}$ dies by passing through a number of
bottlenecks and favorable periods.

\bigskip

\section{Reduced BPREs}

The next stage in studying branching processes is to investigate the
structure of their genealogical trees. For $0\leq k\leq m$ let $Z_{k,m}$ be
the number of particles at moment $k$ in the process $\left\{ Z_n,n\geq
0\right\} $ each of which has a nonempty offspring at moment $m$. The tuple $%
\left\{ Z_{k,m},0\leq k\leq m<\infty \right\} $ is called the reduced
branching process which gives the number of branches in the genealogical
tree of the population modeled by the branching process $\left\{ Z_n,n\geq
0\right\} .$

\subsection{Annealed approach}

The first results for reduced BPREs with linear-fractional PGFs and under
the annealed approach were established by Borovkov and Vatutin \cite{BV} for
the critical case and by Fleischmann and Vatutin \cite{FV} for all three
types of subcritical BPREs. For the critical case Vatitin \cite{V01} proved
under Condition \ref{A1(a)} and some additional conditions that 
\begin{equation}
\mathcal{L}\left( \frac{1}{\sigma \sqrt{n}}\ln Z_{nt,n},\,t\in \lbrack 0,1]\,%
\Big |\,Z_{n}>0\right) \overset{\mathcal{D}}{\rightarrow }\mathcal{L}\left(
\inf_{t\leq u\leq 1}W_{u}^{+},\,t\in \lbrack 0,1]\right) ,  \label{Con1}
\end{equation}%
where $W_{t}^{+},\,t\in \lbrack 0,\infty )$, is the Brownian meander and the
convergence in distribution holds in Skorokhod topology in space $D[0,1]$.
This convergence was established in \cite{BV} for the case of
linear-fractional BPREs under stronger moment assumptions.

Let $B_n=\max \{m<n:Z_{m,n}=1$\}. The difference $d_{n}=n-B_n$ is called the
time to the most recent common ancestor for all individuals existing at time 
$n$ (MRCA$_n$). For ordinary Galton-Watson processes conditioned on
survival, the classical results by Zubkov \cite{Zub} state that in the
critical case $d_{n}$ is asymptotically uniformly distributed over $[0,n]$
while in the subcritical case $d_{n}$ is asymptotically finite. We have a
quite different situation in the random environment setting.

For the annealed approach convergence (\ref{Con1}) states, roughly speaking,
that the number of individuals $Z_{nt,n}$ grows as $\exp \{\sqrt{n}%
\inf_{t\leq u\leq 1}W_{u}^{+}\}$. Recalling that $\mathbf{P}(\inf_{t\leq
u\leq 1}W_{u}^{+}>0)=1$ for any $t\in \left( 0,1\right] ,$ we conclude that
in the annealed setting conditionally on survival the MRCA$_n$ has lived at
early time of order $o\left(n\right) $. In fact, as shown by Borovkov and
Vatutin \cite{BV} for the linear-fractional case, the MRCA$_n$ with positive
probability is the initial individual!

Studying the subcritical BPREs with geometric offspring distributions
Fleischmann and Vatutin \cite{FV} have found that for the intermediately and
strongly subcritical processes the MRCA$_n$ is located not too far from the
point of observation $n$, while for the weakly subcritical case it may be
located either nearby of the point of observation or at the very beginning
of the evolution of the process.

\subsection{Quenched results for the critical reduced BPREs}

This section contains a summary of results established by Vatutin and
Dyakonova in \cite{VD07} and \cite{VD2008} for the critical reduced BPREs
under the quenched approach. The results formulated in this section are
proved under \eqref{dpar} and Conditions \ref{A2}, \ref{A2vd}, if not
explicitely stated otherwise.

\subsubsection{Time to the MRCA and afterwards}

Vatutin and Dyakonova \cite{VD07} proved that in the critical case the MRCA$%
_{n}$ is located "not too far" from the moment of the global minimum $\tau
(n)$ of the ARW on $[0,n]$ in that for any $m\ \in \ \mathbb{Z}$ 
\begin{equation*}
\mathcal{P}\left( B_{n}=\tau (n)+m\right) \overset{d}{\rightarrow }r_{m},\
\sum_{m=-\infty }^{\infty }r_{m}=1,
\end{equation*}%
and with probability 1%
\begin{equation*}
\lim_{m\rightarrow \infty }\lim_{n\rightarrow \infty }\mathcal{P}\left(
|B_{n}-\tau (n)|\leq m\right) =1.
\end{equation*}%
Recalling \eqref{converarc} we clearly see the difference with the ordinary
Galton-Watson critical branching processes where the time to the most recent
common ancestor, scaled by $n$, is uniformly distributed on $[0,1]$.

Also in \cite{VD07} it was proven that the finite-dimensional quenched
distributions of the process $\left\{ Z_{\tau (n)+m,n},m\in \mathbb{Z}%
\right\} $ conditioned on $Z_{n}>0$ weakly converge to the
finite-dimensional distributions of a Galton-Watson branching process
evolving in an inhomogeneous random environment. The limit process starts at 
$-\infty $ by a single individual. The founder of the population dies at a
moment $m\in \left( -\infty ,+\infty \right) $ $\ $\ with probability $r_{m}$
producing at least one offspring in accordance with the PGF $a_{m}(s)$
defined as the weak limit 
\begin{equation*}
\mathcal{E}\left[ s^{Z_{\tau (n)+m+1,n}}|Z_{\tau (n)+m,n}=1\right] \overset{d%
}{\rightarrow }a_{m}(s),s\in \lbrack 0,1),
\end{equation*}%
as $n\rightarrow \infty $. The next generation particles reproduce
independently according to the PGF $a_{m+1}(s).$ And so on. Observe that
this description is in a sharp contrast with the respective limit process
for the ordinary critical reduced Galton-Watson processes obtained in \cite%
{FZ}.


The properties of reduced processes far to the right of $\tau (n)$ were
studied in \cite{VD07} as well. Using the scaling function 
\begin{equation*}
\beta _{n}\left( k\right) :=\frac{1}{\mathcal{E}\left[ Z_{k,n}|\,Z_ n>0%
\right] }, \quad 0\leq k\leq n,  \label{defbeta1}
\end{equation*}
it was shown that as $n\rightarrow \infty$ and $m\rightarrow \infty$ with $%
\tau(n)+m\le n$ 
\begin{equation*}
\mathcal{E}\left[ e^{-\lambda Z_{\tau (n)+m,n}\beta _{n}\left( \tau
(n)+m\right) }\Big|\;Z_n>0\right]\overset{d}{\rightarrow }\Phi (\lambda ),
\end{equation*}%
where $\Phi (\lambda )$ is the same random Laplace transform as in (\ref%
{ddistcont2}). Moreover, if $0\leq m_{1}<m_{2}<...<m_{p}\leq n-\tau (n)$,
then as $m_{1}\rightarrow \infty $ and $n\rightarrow \infty $ 
\begin{equation*}
\mathcal{E} \left[ \prod_{i=1}^{p}e^{-\lambda _{i}Z_{\tau (n)+m_{i},n}\beta
_{n}(\tau (n)+m_{i})}\left\vert Z_n>0\right. \right]\overset{d}{\rightarrow }%
\,\Phi \left( \sum_{i=1}^{p}\lambda _{i}\right) .
\end{equation*}

\subsubsection{Reduced processes at nonrandom times}

Properties of the reduced process $Z_{nt,n}$ for \textit{nonrandom } $t\in
(0,1)$ were also described in \cite{VD07}. It was shown that 
\begin{equation*}
\mathcal{E}\left[ e^{-\lambda Z_{nt,n}\beta _{n}\left( nt\right) }\Big|%
\;Z_n>0\right] \overset{d}{\rightarrow }\Phi \left( \lambda \right)
1_{\left\{ \tau \leq t\right\} }+e^{-\lambda }1_{\left\{ \tau >t\right\}
},\;n\rightarrow \infty ,  \label{DD3}
\end{equation*}%
where $\tau $ is a random point having a Beta$(1-\rho ,\rho )$ distribution
as in \eqref{converarc}  and being independent from $\Phi \left( \lambda
\right) $. Furthermore, for any $\lambda _1\geq 0,\ldots,\lambda_k\ge0$ and $%
0<t_{1}<...<t_{k}\leq 1$, as \ $n\rightarrow \infty $, 
\begin{equation*}
\mathcal{E}\left[ \exp \left\{ -\sum_{i=1}^{k}\lambda _{i}Z_{nt_{i},n}\beta
_{n}\left( nt_{i}\right) \right\} \,|\,Z_{n}>0\right] \overset{d}{%
\rightarrow }\,\Phi \left( \sum_{i=i_{\,\tau }}^{k}\lambda _{i}\right)
e^{-\sum_{i=1}^{i_{\tau }-1}\lambda _{i}},  \label{obsvid}
\end{equation*}%
where $i_{\,\tau }=\min \{i:\,\tau <t_{i}\,\}.$

The corresponding limiting process $\left\{ W_t,0\leq t\leq 1\right\} $
follows a simple pattern (conditioned on a given environmental development).
Until time $\tau$ (corresponding to the time to the MRCA) we have $W_t=1$
and after time $\tau$ we have $W_t\equiv Y$, where the random value $Y$ has $%
\Phi (\lambda )$ as its Laplace transform. Recall that in the
linear-fractional case $Y$ is exponentially distributed with parameter 1.


Concerning the scaling function $\beta _{n}(nt)$ it is known from \cite{VD07}
that under Conditions \ref{B1}, \ref{A2}, \ref{A2vd} for any $t\in (0,1)$ 
\begin{equation*}
\mathcal{L}\left( e^{S_{\tau (nt,n)}-S_{\tau (n)}}\beta _{n}(nt)|\tau
(n)<nt\right) \overset{d}{\rightarrow }\mathcal{L}\left( \zeta ^{\ast
}\right), \quad n\rightarrow \infty ,
\end{equation*}%
where $\zeta ^{\ast }$ is a proper random variable. It follows that if the
global minimum time $\tau (n)$ of the ARW occurs prior to $nt$ then $Z_{nt,n}
$ is of order $e^{S_{\tau (nt,n)}-S_{\tau (n)}}$ (up to a random multiplier
separated from zero and infinity).

\subsubsection{Reduced BPREs as random fields}

Vatutin and Dyakonova \cite{VD2008} obtained quenched limit theorems for the
random fields $\{Z_{nt_{1},nt_{2}}\}_{0<t_{1}<t_{2}<1}$ as $Z_{n}>0$ and $%
n\rightarrow \infty $. The answers essentially depend on the position of $%
\tau (n)$ relative to the interval $[nt_{1},nt_{2}]$: as $n\rightarrow
\infty $ 
\begin{align}
& \mathcal{L}\left( \mathcal{E}[e^{-\lambda Z_{nt_{1},nt_{2}}\beta
_{n}(t_{1},t_{2})}|Z_{n}>0]\,\Big |\,\tau (n)<nt_{1}\right) \overset{d}{%
\rightarrow }\Phi \left( \lambda \right) ,  \label{K1} \\
& \mathcal{L}\left( \mathcal{E}[e^{-\lambda Z_{nt_{1},nt_{2}}\beta
_{n}(t_{1},t_{2})}|Z_{n}>0]\,\Big |\,\tau (n)>nt_{2},\tau (nt_{2})\leq
nt_{1}\right) \overset{d}{\rightarrow }\Psi \left( \lambda \right) ,
\label{K2} \\
& \mathcal{L}\left( \mathcal{E}[s^{Z_{nt_{1},nt_{2}}}|Z_{n}>0]\Big |%
\,nt_{1}\leq \tau (n)\leq nt_{2}\right) \overset{d}{\rightarrow }s,
\label{K3} \\
& \mathcal{L}\left( \mathcal{E}[s^{Z_{nt_{1},nt_{2}}}|Z_{n}>0]\Big |\,\tau
(n)>nt_{2},\tau (nt_{2})>nt_{1}\right) \overset{d}{\rightarrow }s,
\label{K4}
\end{align}%
where $\Phi \left( \lambda \right) $ is from (\ref{ddistcont2}), $\Psi
\left( \lambda \right) $ is from (\ref{addit}), and 
\begin{equation*}
\beta _{n}(t_{1},t_{2}):=\frac{1}{\mathcal{E}\left[ Z_{nt_{1},nt_{2}}|%
\,Z_{nt_{2}}>0\right] }=\frac{1-f_{0,nt_{2}}(0)}{e^{S_{nt_{1}}}\left(
1-f_{nt_{1},nt_{2}}(0)\right) }.
\end{equation*}

According to Theorem 1 in \cite{VD} 
\begin{equation*}
\beta _{n}(t_{1},t_{2})e^{S_{\tau (nt_{1},nt_{2})}-S_{\tau (nt_{1})}}\overset%
{d}{\rightarrow }\hat{\Theta},\ n\rightarrow \infty ,
\end{equation*}%
where the random variable $\hat{\Theta},$ is positive and finite with
probability 1. Observe that in the cases \eqref{K3} and \eqref{K4} we have $%
S_{\tau (nt_{1},nt_{2})}\leq S_{\tau (nt_{1})}$. Results \eqref{K3} and %
\eqref{K4} say that in these cases all individuals existing at moment $nt_{2}
$ are descendants from a single individual existing at moment $nt_{1}$.
Relations \eqref{K1} and \eqref{K2} describe the cases when the
corresponding number of ancestors is large with $Z_{nt_{1},nt_{2}}\asymp
e^{S_{\tau (nt_{1},nt_{2})}-S_{\tau (nt_{1})}}$.

\end{document}